%% file: pulse3.tex
\documentclass[10pt,a4paper,leqno]{amsart}

\usepackage{amsmath,amssymb,amsthm,epsfig,euscript}

\def\eps
{\varepsilon}

\def\R{{\mathbb R}}

\def\virgp{\raise 2pt\hbox{,}}

\def\bu{{\bf u}}

\def\ut{\underline t}
\def\ur{\underline r}

\def\Eq#1#2{\mathop{\sim}\limits_{#1\rightarrow#2}}
\def\Tend#1#2{\mathop{\longrightarrow}\limits_{#1\rightarrow#2}}

\def\d{{\partial}}
\def\e{\varepsilon}

\theoremstyle{plain}
\newtheorem{theo}{Theorem}[section]
\newtheorem{lem}[theo]{Lemma}
\newtheorem{cor}[theo]{Corollary}
\newtheorem{prop}[theo]{Proposition}

\theoremstyle{definition}
\newtheorem{defin}[theo]{Definition}

\theoremstyle{remark}
\newtheorem*{rema}{Remark}

\newtheorem{ex}[theo]{Example}

\numberwithin{equation}{section}

\begin{document}

\title[Focusing of Spherical Nonlinear
Pulses in ${\mathbb R}^{1+3}$, III]{Focusing of Spherical Nonlinear
Pulses in ${\mathbb R}^{1+3}$, III. Sub and Supercritical cases}  
\author[R. Carles]{R{\'e}mi Carles}
\address{MAB, UMR CNRS 5466\\
Universit{\'e} Bordeaux 1\\ 351 cours de la Lib{\'e}ration\\ 33~405 Talence
cedex\\ France}
\email{carles@math.u-bordeaux.fr}
\author[J. Rauch]{Jeffrey Rauch}
\address{Department of
Mathematics\\ University of Michigan\\ Ann Arbor\\ MI 48109\\ USA}
\email{rauch@math.lsa.umich.edu}  
\begin{abstract}
We study the validity of geometric optics in $L^\infty$
for nonlinear wave
equations in three space dimensions whose solutions, pulse like, focus at a
point. If the amplitude of the initial data is
subcritical, then no nonlinear effect occurs at leading order. 
If the amplitude of the initial data is sufficiently big, 
strong nonlinear effects occur; we study the cases where the equation
is either dissipative or accretive. 
When the equation is dissipative,
pulses are absorbed before reaching the focal point. When the equation
is accretive,  the family of 
 pulses becomes unbounded. 
\end{abstract}
\maketitle


\input{intro}

\input{approx}

\input{exist}

\input{sub}

\input{super}

\bibliographystyle{amsalpha}
\bibliography{../carles}

\end{document}

%% file: intro.tex
\section{Introduction}
\label{sec:intro}

This paper is the last of a series of three, after \cite{CR} and
\cite{CR2}. In these three papers, we consider the asymptotic behavior
as $\e\to 0$ of solutions 
of the initial value problem
\begin{equation}
\label{eq:pbinit}
\left\{
  \begin{split}
   \square \bu^\e + a |\d_t \bu^\e |^{p-1}\d_t \bu^\e&\ =\ 0, \qquad 
(t,x)\in [0,T]\times\R^3, \\
   \bu^\e\big|_{ t=0}&\ =\ 
        \e^{J+1} U_0\left(r\,,\frac{r-r_0}{\e}\right),\\ 
   \d_t \bu^\e\big|_{ t=0}&\ =\  
        \e^J U_1\left(r\,,\frac{r-r_0}{\e}\right),
  \end{split}
\right.
\end{equation}

\noindent
where $\square := \d_t^2-\Delta_x$, 
 $a$ is a complex number, $r=|x|$, $r_0 >0$, and,
$1< p<\infty$.
The functions $U_0$ and $U_1$
are real-valued, infinitely differentiable, bounded, and 
there is a $z_0>0$ so that
 for all $r\geq 0$,
\begin{equation}
\label{eq:Usupport}
\operatorname{supp}U_j(r,.)\subset [-z_0,z_0].
\end{equation}
The last assumption implies that at time $t=0$ the solutions are
families of spherical pulses
 supported in a 
$O(\e)$ neighborhood of  $r=r_0$. The initial data are spherically
symmetric, thus in the limit $\e \to 0$, a caustic is formed, reduced
to the focal point $(t,x)=(r_0,0)$. 
Before
going further into details, we rescale our parameters as in
\cite{CR2}. Introduce
$\eps^{-J}\bu^\eps=:u^\eps$ instead of $\bu^\eps$ so
that the solutions have 
derivatives of order $O(1)$ away from the caustic.
Define $\alpha \ :=\ (p-1)J$. The initial
value problem \eqref{eq:pbinit} is transformed to

\begin{equation}
\label{eq:pb}
\left\{
  \begin{split}
   &\square u^\eps + a \,\eps^\alpha |\d_t u^\eps |^{p-1}\,\d_t u^\eps=0, \ \ 
(t,x)\in [0,T]\times\R^3, \\
   &u^\eps\big|_{ t=0}= \eps U_0\left
   ( r,\frac{r-r_0}{\eps}\right),
\quad 
\quad
   \d_t u^\eps\big|_{ t=0}=   U_1\left( r,\frac{r-r_0}{\eps}\right)\,.
  \end{split}
\right.
\end{equation}
In \cite{CR2}, formal arguments,
inspired by the linear case $a=0$, led to the following distinctions,
in the spirit of those computed formally in \cite{HK87},
\bigskip
\begin{center}
\begin{tabular}[c]{|l|c|c|c|}
\hline
 &$\alpha+2>p$  & $\alpha+2=p$ & $\alpha+2<p$ \\
\hline
$\alpha >0$ &linear caustic, & nonlinear caustic,& supercritical caustic,\\
&linear propagation   & linear propagation & linear propagation \\
\hline
$\alpha =0$&linear caustic, &nonlinear caustic,& supercritical caustic,\\
&nonlinear propagation   & nonlinear propagation & nonlinear propagation\\
\hline
\end{tabular}
\end{center}

\medskip

In \cite{CR}, we studied the case ``linear caustic, nonlinear
propagation''; we proved that nonlinear geometric optics provides a
good approximation of $\d_t \bu^\e$ away from the focal point
$(t,r)=(r_0,0)$, and that the nonlinear term is negligible near the
focus. In \cite{CR2}, we analyzed the case ``nonlinear caustic, linear
propagation''. In some sense, it is the exact opposite of the previous
case;  the nonlinear term is negligible outside the focal point, but
has a relevant influence near $(r_0,0)$, which is described by a
nonlinear scattering operator. Moreover, this scattering operator
broadens the pulses (at least if $U_0$ and $U_1$ are small), which
leave the focus with algebraically decaying tails. 

In this paper, we discuss the remaining cases of the above table. In
the last three cases, we treat only the case of $a$ real, that is when 
Eq.~\eqref{eq:pb} is dissipative or accretive (in particular, we do not 
treat the case of
conservative equations).  
\bigskip

The first case, ``linear caustic, linear propagation'', suggests that
the nonlinear term is everywhere negligible; we prove that this is
so. In the last three cases, we assume that the coupling constant $a$
is real. For the ``supercritical caustic'' cases, strong nonlinear effects
are expected near the focal point. We prove that when
the equation is dissipative ($a>0$), then the
dissipation is so strong near the focus that the pulses are
absorbed. This result is the pulse analogue of \cite{JMRTAMS95} and
\cite{JMRMemoir}, which proved absorption in the case of wave trains
 (the initial profiles $U_j$ are assumed to be periodic with
respect to their last variable instead of compactly supported).  
More precisely, in \cite{JMRTAMS95} and \cite{JMRMemoir}, it is proved
that the exact solution $u^\eps$ of 
the dissipative ($a>0$) wave equation \eqref{eq:pbinit} with $J=0$,
is approximated as follows,
$$\d_t u^\eps(t,x)\Eq \eps 0
\d_t\underline{u}(t,x)+U_-\left(t,x,\frac{t+|x|}{\eps} \right) +
U_+\left(t,x,\frac{t-|x|}{\eps} \right) ,$$
where the profiles  $U_\pm$ are periodic with respect to their last
variable, with mean value zero. The absorption of oscillations is
given by $U_\pm \equiv 0$ past the caustic. Thus, only
the average term remains, included in $\d_t\underline{u}$.
For an almost periodic function, the notion of average is given by
$$\underline{f}=\lim_{T\rightarrow +\infty}\frac{1}{2T}\int_{-T}^T
f(\theta)d\theta. $$
When  $f$ is compactly supported, the above limit is zero, and pulses
formally have mean value zero. Thus, the absorption of pulses is the
formal analogue of the absorption of oscillations. 

Our
present framework makes it possible to analyze very precisely the
corresponding phenomenon for pulses; they are absorbed when
approaching the caustic, that is even \emph{before reaching it}. We 
prove that this phenomenon occurs for the last three cases of the
table, ``supercritical caustic, linear/nonlinear propagation'' and
``nonlinear caustic, nonlinear propagation''. 

The present paper along with \cite{CR} and \cite{CR2} prove
 that the distinctions derived formally in \cite{CR2} and
recalled in the above tables are correct. Let us give an
interpretation of these results when the nonlinearity is fixed, and
when one modulates the amplitude of the initial data in
\eqref{eq:pbinit}. Consider a fixed $p>2$, and modify the value of
$J$. For a unified, complete, presentation, we assume that the equation 
is dissipative,
$a>0$.  
\begin{itemize}
\item If $J>\frac{p-2}{p-1}$, then the pulse is not affected by the
nonlinearity at leading order. It remains too small to ignite the
nonlinearity. 
\item If $J=\frac{p-2}{p-1}$, then  the
nonlinear term is negligible away from the focus, but the caustic
crossing, described by a scattering operator, has enlarged the support
of the pulse, and decreased its 
amplitude. The pulse is too small to see the nonlinearity outside the
focal point, but the amplification near the caustic makes the
nonlinear term relevant there. 
\item If $J<\frac{p-2}{p-1}$, then  the pulse
is absorbed at the focus. It is sufficiently big to make the nonlinear
effects so strong that the dissipation is complete before the focus.
\end{itemize}
When $p=2$, the nonlinear term is negligible if $J>0$, and if $J=0$,
the pulses are absorbed. If $1<p<2$, the same method would prove that
the nonlinear term is 
negligible if $J>0$, and the case $J=0$ was treated in \cite{CR}.\\

Before stating precisely our results, we make a change of
unknown, as in \cite{CR} and \cite{CR2}. Since the
initial data are spherical, so is the solution.
With the usual
abuse of notation,
$$
u^\eps(t,x) \ = \ u^\eps(t,|x|),
\qquad
u^\eps(t,|x|)\ \in \  C^\infty_{{\rm even\ in} \ r}(\R_t\times\R_r)\,.
$$
 Introduce $v^\eps:=(v^\eps_-,v^\eps_+)$ where
\begin{equation}
\label{eq:related}
\tilde u^\eps(t,r) := ru^\eps(t,r),
\qquad v_\mp^\eps := (\d_t \pm\d_r)\tilde u^\eps\,,
\qquad
v_\mp^\eps\in C^\infty (\R_t\times\R_r)
\,.
\end{equation}
Then \eqref{eq:pb} becomes

\begin{equation}
\label{eq:pbreduit}\left\{
  \begin{split}
      &(\d_t\pm \d_r)v_\pm^\eps = \eps^{\alpha}r^{1-p}g(v_-^\eps +v_+^\eps),
\qquad
g(y):=-a2^{-p} |y|^{p-1}y\,,
\\
 &  (v_-^\eps +v^\eps_+)\big|_{r=0}=0\,,\\
&v_\mp^\eps\big|_{t=0} = P_\mp\left( r,\frac{r-r_0}{\eps}\right) \pm \eps
P_1\left( r,\frac{r-r_0}{\eps}\right)\,,
\end{split}\right.
\end{equation}
where
\begin{equation*}
\begin{split}
 P_\mp(r,z)&
\ :=\ 
rU_1(r,z)\pm r\d_z U_0(r,z)\,, \\
 P_1(r,z)&
\ :=\ 
U_0(r,z)+r\d_rU_0(r,z)\,.
\end{split}
\end{equation*}
We prove asymptotics for $v^\eps$; asymptotics for $\d_t u^\eps$
are deduced by \eqref{eq:related}, 
$$\d_t u^\eps (t,r) =\frac{v_-^\eps + v^\eps_+}{2r}.$$
For the sub-critical case, introduce the solution of the linear
equation
\begin{equation}
\label{eq:libre}\left\{
  \begin{split}
      (\d_t\pm \d_r)(v_\pm^\eps)_{\rm free} &= 0,\\
   \left((v_-^\eps)_{\rm free} +(v^\eps_+)_{\rm
   free}\right)\big|_{r=0}&=0\,,\\
  (v_\mp^\eps)_{\rm free}\big|_{t=0} &= P_\mp
\left( r,\frac{r-r_0}{\eps}\right) \,.
\end{split}\right.
\end{equation}
It is given explicitly by the formulae,
\begin{equation*}
\begin{split}
(v_-^\eps)_{\rm free}(t,r)&=P_-\left(r+t,\frac{r+t-r_0}{\e}\right),\\
(v_+^\eps)_{\rm free}(t,r)&=P_+\left(r-t,\frac{r-t-r_0}{\e}\right)-
P_-\left(t-r,\frac{t-r-r_0}{\e}\right).
\end{split}
\end{equation*}
The pulse $(v_-^\eps)_{\rm free}$ corresponds to an incoming wave, and
$(v_+^\eps)_{\rm free}$ is the sum of two outgoing waves, one from
$P_+$, and the other from the focusing of the incoming wave.  
\begin{theo}[Sub-critical case]
\label{th:sub}
Assume that $\alpha > \max (0,p-2)$. \\
$\bullet$ If $\alpha>p-2>0$, then there exists $\e_0>0$ such that for
any $\e\in]0,\e_0]$,
\eqref{eq:pbreduit} has a unique, global, solution $v^\eps\in
C^1([0,\infty[\times\R_+)$. Moreover, the following asymptotics holds
in $L^\infty(\R_+\times\R_+)$, 
\begin{equation*}
\begin{split}
v^\e_\pm (t,r)& = (v_\pm^\eps)_{\rm free}(t,r) +
O\left(\eps^{\min(1,\alpha+2-p)}\right),\\ 
\e\d_t v^\e_\pm (t,r)& = \e\d_t(v_\pm^\eps)_{\rm free}(t,r) +
O\left(\eps^{\min(1,\alpha+2-p)} 
\right).
\end{split}
\end{equation*}
$\bullet$ If $\alpha>0$ and $1<p\leq 2$, let $T>0$. Then there exists
$\e(T)>0$ such that for 
any $\e\in]0,\e(T)]$,
\eqref{eq:pbreduit} has a unique solution $v^\eps\in
C^1([0,T]\times\R_+)$. Moreover, the following asymptotics holds
in $L^\infty([0,T]\times\R_+)$, 
\begin{equation*}
\begin{split}
v^\e_\pm (t,r)& = (v_\pm^\eps)_{\rm free}(t,r) +
O\left(\eps^{\min(1,\alpha)}\right) \ \ \left(
O\left(\eps+ \e^\alpha |\log \e|\right) \textrm{ if }p=2\right),\\ 
\e\d_t v^\e_\pm (t,r)& = \e\d_t(v_\pm^\eps)_{\rm free}(t,r) +
O\left(\eps^{\min(1,\alpha)}\right) \ \ \left(
O\left(\eps+ \e^\alpha |\log \e|\right) \textrm{ if }p=2\right).
\end{split}
\end{equation*}
\end{theo}

\begin{theo}[Super-critical case]
\label{th:super}
Assume that $0\leq \alpha <p-2$ or $\alpha =0=p-2$. \\
$\bullet$ If the equation is 
dissipative, $a >0$, the pulses are absorbed before
reaching the focus. If $T\geq r_0$,
\begin{equation*}
\limsup_{\eps \rightarrow 0}\left(\|v_-^\eps(T)\|_{L^\infty(0\leq r
    \leq T)}
+\|v_+^\eps(T)\|_{L^\infty(0\leq r
    \leq T)}\right) =0.
\end{equation*}
More precisely, for $\lambda >0$, define $T(\lambda,\e)$ as follows;
if $0\leq \alpha 
<p-2$, then $T(\lambda,\e):=r_0 -z_0\e -\lambda \e^{\alpha/(p-2)}$,
and if $\alpha =0=p-2$, then  $T(\lambda,\e):=r_0 -z_0\e -\lambda$. 
For any $T=T(\e)\geq T(\lambda,\e)$, 
\begin{equation*}
\lim_{\lambda \to 0}
\limsup_{\eps \rightarrow 0}\left(\|v_-^\eps(T)\|_{L^\infty(0\leq r
    \leq T)}
+\|v_+^\eps(T)\|_{L^\infty(0\leq r
    \leq T)}\right) =0.
\end{equation*}
$\bullet$ If the equation is accretive, $a<0$, there exists $T^*\leq
r_0$ such that 
the family $(v_-^\eps, v^\e_+)$ is not bounded in
$L^\infty([0,T^*]\times\R_+)^2$. 
\end{theo}
\begin{rema}
For $0<\e \leq 1$ and $\lambda$ positive,
$T(\lambda,\e)<r_0-z_0\e$, so the first part of the theorem shows that the
absorption mechanism takes place before the incoming wave reaches the
focus. Indeed, the pulses are initially supported in $\{|r-r_0|\leq
z_0\e\}$, so by finite speed of propagation, they do not reach the
origin before $t=r_0-z_0\e$ (see Fig.~\ref{fig:ogg}). 
\begin{figure}[htbp]
\begin{center}
\input{ogg.pstex_t}
\caption{Geometry of the propagation in the super-critical case.}
\label{fig:ogg}
\end{center}
\end{figure}
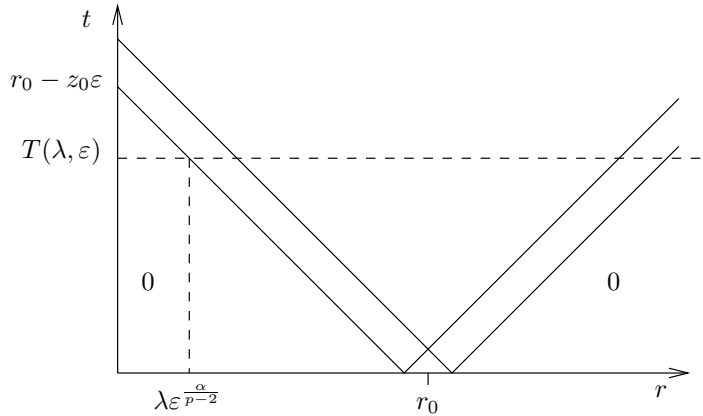 
Notice that in
the dissipative case, our estimates are in $\{t\geq r\}$, which
includes the  focusing region, and its domain of influence. The key to
the proof of Theorem~\ref{th:super} is to construct and approximate
solution which is more accurate than the approximation of nonlinear
geometric optics and which permits us to penetrate with high accuracy
to small distances from the focal point $r=0$.
\end{rema}

\begin{rema}
If we considered wave trains, that is
$U_j(r,.)$ periodic instead of compactly supported,
Theorem~\ref{th:sub} would still hold. The
proof we 
give works in both cases. On the other hand, the proof of
Theorem~\ref{th:super} relies on the compact support
assumption. We construct approximate solutions that solve ordinary
differential equations along the rays of geometrical optics (see
Sect.~\ref{sec:super}), and can be computed explicitly. In the
dissipative case, these approximate solutions are absorbed before they
stop being good approximations, proving thereby the absorption of the
exact solution. In the case of wave trains, the computation of the
counterpart of these approximations is a project for the future.  
\end{rema}

\begin{rema}
As we mentioned above, in the supercritical 
cases, our framework is restricted. We assume that the coupling 
constant $a$ is real, while we do not
make this assumption in Theorem~\ref{th:sub}. It would be interesting to know
what happens when, for instance, $a$ is pure imaginary; no absorption can 
happen, for the equation in that case is conservative. A partial answer is 
given in \cite{CL}, in the case $\alpha =0=p-2$, on a system which is a 
simplified model for \eqref{eq:pb}: an arbitrary  phase shift appears, 
varying like $\log \eps$. In the supercritical framework, $\alpha =0$ and 
$p>2$, one may expect even more pathological behaviors. 
\end{rema}

In Section~\ref{sec:approx}, we prove two stability results. 
In Section~\ref{sec:exist}, we discuss existence
results, and establish estimates for the sub-critical
case. Theorem~\ref{th:sub} is proved in Section~\ref{sec:sub}, and
Theorem~\ref{th:super} is proved in Section~\ref{sec:super}.

The results of Theorem~\ref{th:super} were announced in
\cite{CRcras1}.

%% file: ogg.pstex_t
\begin{picture}(0,0)%
\includegraphics{ogg.pstex}%
\end{picture}%
\setlength{\unitlength}{3947sp}%
\begingroup\makeatletter\ifx\SetFigFont\undefined%
\gdef\SetFigFont#1#2#3#4#5{%
  \reset@font\fontsize{#1}{#2pt}%
  \fontfamily{#3}\fontseries{#4}\fontshape{#5}%
  \selectfont}%
\fi\endgroup%
\begin{picture}(4362,2611)(1,-1835)
\put(  1,239){\makebox(0,0)[lb]{\smash{\SetFigFont{10}{12.0}{\rmdefault}{\mddefault}{\updefault}$r_0-z_0\e$}}}
\put( 76,-211){\makebox(0,0)[lb]{\smash{\SetFigFont{10}{12.0}{\rmdefault}{\mddefault}{\updefault}$T(\lambda,\e)$}}}
\put(3751,-1036){\makebox(0,0)[lb]{\smash{\SetFigFont{10}{12.0}{\rmdefault}{\mddefault}{\updefault}$0$}}}
\put(451,614){\makebox(0,0)[lb]{\smash{\SetFigFont{10}{12.0}{\rmdefault}{\mddefault}{\updefault}$t$}}}
\put(4051,-1711){\makebox(0,0)[lb]{\smash{\SetFigFont{10}{12.0}{\rmdefault}{\mddefault}{\updefault}$r$}}}
\put(2551,-1786){\makebox(0,0)[lb]{\smash{\SetFigFont{10}{12.0}{\rmdefault}{\mddefault}{\updefault}$r_0$}}}
\put(901,-1786){\makebox(0,0)[lb]{\smash{\SetFigFont{10}{12.0}{\rmdefault}{\mddefault}{\updefault}$\lambda\eps^{\frac{\alpha}{p-2}}$}}}
\put(826,-1036){\makebox(0,0)[lb]{\smash{\SetFigFont{10}{12.0}{\rmdefault}{\mddefault}{\updefault}$0$}}}
\end{picture}

%% file: approx.tex
\section{General stability results}
\label{sec:approx}
In this section, we state two general approximation arguments. The
first one  will
allow us to prove Theorem~\ref{th:sub} and the second one,
Theorem~\ref{th:super}. Our first result is an easy estimate, proved
in \cite{CR2}.
\begin{defin}
For $t>0$, we denote by $\Gamma_-^t$ (resp. $\Gamma_+^t$) the set of
all speed minus one (resp. plus one) characteristics connecting points
on the initial line $\{t=0\}$ to  points at time $t$. We also denote
$\Gamma^t = \Gamma_-^t \cup \Gamma_+^t$.
\end{defin}
For a characteristic $\gamma \in \Gamma^t$, we use the convention that 
$$\int_\gamma f$$
stands for the integral of $f$ along $\gamma$, parameterized by the
time variable,
$$\int_\gamma f= \int_0^t f(s,r(s))ds,$$
where $(s,r(s))_{s\in [0,t]}$ is a parametrization of $\gamma$. In
particular, if $f$ is nonnegative, then so is the above integral. 
\begin{lem}[\cite{CR2}, Lemma~2.1]  
\label{lem:linearest}
Suppose that $w$ and $f=(f_+ , f_-)$
are bounded continuous  functions on
$[0,T]\times [0,\infty[$ satisfying in the sense of distributions
\begin{equation*}
(\d_t\pm\d_r)w_\pm = f_\pm\,,
\qquad w_+(t,0) + w_-(t,0) = 0\,,\quad
{\rm for}\ \ 0\le t\le T\,.
\end{equation*}
Denote by 
$$
M_\pm(t) \ :=\ 
\|w_\pm(t)\|_{L^\infty([0,\infty[)}\,.
$$
Then for $0\le t\le T$ one has
\begin{equation*}
M_\pm(t)
\ \le \ 
\max \{M_+(0), M_-(0)\}
\ +\
\sup_{\gamma_- \in \Gamma_-^t}
\int_{\gamma_-}|f_-|+\
\sup_{\gamma_+ \in \Gamma_+^t}
\int_{\gamma_+}|f_+|.
\end{equation*}
\end{lem}
To prove Theorem~\ref{th:super}, we use a result in
the spirit of Gronwall's lemma. The assumption $0<T<\delta$ guarantees
that the support of the solution does not touch $\{ r=0 \}$. 
\begin{prop}
\label{prop:approx}
Suppose that $0<T<\delta$, and $w=(w_-,w_+)\in C\cap L^\infty
([0,T]\times \R_+)$ satisfies
\begin{equation}
\label{eq:ref}
\left\{
\begin{split}
(\d_t \pm \d_r)w_\pm & =f_\pm (t,r)(w_- + w_+) + S_\pm,\\
w_{\pm \mid t=0}&= w_{0\pm},
\end{split}
\right.
\end{equation}
with $ \operatorname{supp}
w_{0\pm} \subset [\delta,+\infty[$. Suppose 
that 
\begin{equation}\label{eq:assum}
C_1=C_1(f):= \int_0^{T}\sup_{\gamma_-\in
\Gamma_-^t} |f_-|dt +\int_0^{T} \sup_{\gamma_+\in
\Gamma_+^t}|f_+|dt \, < \, \infty.
\end{equation}
Then 
\begin{equation*}
\sup_{0\leq t \leq T}\|w_\pm(t)\|_{L^\infty} \leq C_2\sum_{\pm}\Big(
\|w_{0\pm}\|_{L^\infty} + 
\sup_{\gamma_\pm\in
\Gamma_\pm^{T}}\int_{\gamma_\pm} |S_\pm| \Big), 
\end{equation*}
with $C_2=\max ( C_1 e^{2C_1}, C_1^2 e^{3C_1})$.
\end{prop}
\noindent {\bf Warning.} Note that in hypothesis \eqref{eq:assum}, the
supremum is inside the integral. The estimate would not be true with
the supremum outside. 
\begin{proof}
Let $(t,r)\in [0,T]\times ]0,\infty[$, and denote $\gamma_-
=\gamma_-(t,r)$ the characteristic from $(0,t+r)$ to
$(t,r)$. Duhamel's principle for $w_-$ reads
$$w_-(t,r)=w_{0-}(t+r) +\int_{\gamma_-}f_-\times (w_- +
w_+) + \int_{\gamma_-}S_- .$$
Gronwall lemma, along with assumption~\eqref{eq:assum}, yields, for any
$t\in [0,T]$,
\begin{equation}\label{eq:estim1}
\begin{split}
\|w_-(t)\|_{L^\infty} \leq & e^{C_1}\left( 
\|w_{0-}\|_{L^\infty} + \sup_{\gamma_-\in \Gamma_-^t}\int_{\gamma_-}
|f_-. w_+| + \sup_{\gamma_-\in
 \Gamma_-^t}\int_{\gamma_-}|S_-|\right)\\
\leq &e^{C_1}\left(\|w_{0-}\|_{L^\infty} + C_1 \sup_{0\leq s \leq
 t}\|w_+(s)\|_{L^\infty} + \sup_{\gamma_-\in
 \Gamma_-^t}\int_{\gamma_-}|S_-|\right).
\end{split}
\end{equation}
Similarly, 
\begin{equation}\label{eq:*}
\|w_+(t)\|_{L^\infty} \leq  e^{C_1}\left( 
\|w_{0+}\|_{L^\infty} + \sup_{\gamma_+\in \Gamma_+^t}\int_{\gamma_+}
|f_+. w_-| + \sup_{\gamma_+\in
 \Gamma_+^t}\int_{\gamma_+}|S_+|\right).
\end{equation}
For the $w_-$ integrals on the right in \eqref{eq:*}, use estimate
\eqref{eq:estim1}  
to find
$$ \int_{\gamma_\pm}|f_\pm w_-| 
\ \le \
\int_{\gamma_\pm}
|f_\pm|\
e^{C_1}\left(\|w_-(0)\| +
C_1\max_{0\le \tau\le s}\|w_+(\tau)\|+\sup_{\widetilde\gamma_-\in
\Gamma_-^s}\int_{\widetilde\gamma_-} |S_-|\right)ds,
$$
where $\gamma_\pm$ is parameterized by the time $s\in [0,t]$.
Introduce
$$
{\bf m_\pm}(t) :=\sup_{0\le s\le t}\|w_\pm(s)\|\,,
$$
and take the supremum
on $r$ in \eqref{eq:*}
to obtain
\begin{equation*}
\begin{split}
\|w_+(t)\|
\ \le \
C_1&e^{C_1}\Bigg(\sum_{\pm}\left(\|w_\pm(0)\|
+  \sup_{\gamma_\pm\in
\Gamma_\pm^t}\int_{\gamma_\pm} |S_\pm|\right)+\\ 
&+ \sup_{\gamma_+\in
\Gamma^t_+ } \int_{\gamma_+} |f_+(s,r)|\ {\bf m}_+(s) ds
+ \sup_{\gamma_-\in \Gamma^t_-
} \int_{\gamma_-} |f_-(s,r)|\ {\bf m}_+(s) ds\Bigg).
\end{split}
\end{equation*}
Therefore,
\begin{equation*}
\begin{split}
\|{\bf m}_+(t)\|
\ \le \
C_1&e^{C_1}\Bigg(\sum_{\pm}\left(\|w_\pm(0)\|
+  \sup_{\gamma_\pm\in
\Gamma_\pm^t}\int_{\gamma_\pm} |S_\pm|\right)+\\ 
&+ \sup_{\gamma_+\in
\Gamma^t_+ } \int_{\gamma_+} |f_+(s,r)|\ {\bf m}_+(s) ds
+ \sup_{\gamma_-\in \Gamma^t_-
} \int_{\gamma_-} |f_-(s,r)|\ {\bf m}_+(s) ds\Bigg)\\
\ \le \
C_1&e^{C_1}\Bigg(\sum_{\pm}\left(\|w_\pm(0)\|
+  \sup_{\gamma_\pm\in
\Gamma_\pm^t}\int_{\gamma_\pm} |S_\pm|\right)+\\ 
&+ \int_0^t \sup_{\gamma_+\in
\Gamma^t_+ } |f_+(s,r)|\ {\bf m}_+(s) ds
+  \int_0^t \sup_{\gamma_-\in \Gamma^t_-
}|f_-(s,r)|\ {\bf m}_+(s) ds\Bigg).
\end{split}
\end{equation*}
Applying Gronwall's lemma, using assumption \eqref{eq:assum}, yields
\begin{equation*}
{\bf m}_+(T) \leq C_1e^{2C_1}\left(\|w_+(0)\|
+   \|w_-(0)\|+ 
\sup_{\gamma_-\in
\Gamma_-^{T}}\int_{\gamma_-} |S_-| +\sup_{\gamma_+\in
\Gamma_+^{T}}\int_{\gamma_+} |S_+|\right).
\end{equation*}
This inequality, along with \eqref{eq:estim1}, proves the proposition.
\end{proof}
The following example contains the core of the proof of
Th.~\ref{th:super}. 
\begin{ex}\label{ex:gener}
Consider the case 
$$f^\e_\pm(t,r)= \e^\alpha r^{1-p}|\widetilde v^\e_-(t,r)|^{p-1},$$
where $\widetilde v^\e_-$ solves an initial value problem of the form
(compare with \eqref{eq:app} below), 
$$(\d_t -\d_r)\widetilde v^\e_- = F^\e\left( \widetilde v^\e_-
\right)\ ,\ \ \widetilde v^\e_{-\mid
t=0} = P_-\left( 
r,\frac{r-r_0}{\eps}\right).$$
Recall that $\operatorname{supp}P_-(r,.)\subset [-z_0,z_0]$. For 
$w^\e_{0_\pm}$, consider pulses with the same support as $\widetilde 
v^\e_{-\mid t=0}$. Then by finite speed
of propagation, we can take $\delta^\e =r_0-z_0\e$ in
Prop.~\ref{prop:approx}. For $0<\ut<\delta^\e$, the maximum of
$f^\e_\pm$ at time $\ut$ is estimated by
$$\e^\alpha {\ur}^{1-p}\|\widetilde v^\e_-\|^{p-1}_{L^\infty(0\leq t\leq
\ut)},$$ 
where $\ur$ is such that $\ur+\ut=\delta^\e =r_0-z_0\e$. Therefore we
have
\begin{equation*}
\begin{split}
\int_0^{\ut}\sup_{\gamma_-\in
\Gamma_-^t} |f_-^\e|dt +\int_0^{\ut} \sup_{\gamma_+\in
\Gamma_+^t}|f_+^\e|dt &\leq 2\|
\widetilde v^\e_-\|^{p-1}_{L^\infty(0\leq t\leq \ut)} 
\int_0^{\ut}\e^\alpha (\delta^\e -t)^{1-p}dt \\
&\leq C_p\e^\alpha\|
\widetilde v^\e_-\|^{p-1}_{L^\infty(0\leq t\leq \ut)}\times\left\{
  \begin{split}
   (\delta^\e -\ut)^{2-p}& \textrm{ if }p>2,\\
   |\log (\delta^\e -\ut)|& \textrm{ if }p=2.
  \end{split}\right.
\end{split}
\end{equation*}
Thus, when $\widetilde v^\e_-$ remains bounded and $T^\e$ is chosen
so that $\e^\alpha (\delta^\e -T^\e)^{2-p}$ (resp. $\e^\alpha|\log
(\delta^\e -T^\e)| $ if $p=2$) is bounded independent of $\e$, then we
can use Prop.~\ref{prop:approx}. This is the case in particular if
$T^\e=T_{\lambda,\e}$ defined in Th.~\ref{th:super}. 
\end{ex}

%% file: exist.tex
\section{Existence results}
\label{sec:exist}
In this section, we prove two kinds of results concerning the existence of
solutions to \eqref{eq:pbreduit}; local existence in $L^\infty$
before the wave meets the boundary $\{r=0\}$, and global existence in
$W^{1,\infty}$ when the boundary condition has to be taken into
account. The reason appears in Th.~\ref{th:sub} and \ref{th:super}; in
Th.~\ref{th:super}, the phenomena we want to prove occur \emph{before}
the pulses reach the boundary, while Th.~\ref{th:sub} includes
the caustic crossing.

In the first case, we are interested in a problem
\begin{equation}
\label{eq:pbreduit2}\left\{
  \begin{split}
      &(\d_t\pm \d_r)v_\pm^\eps = \eps^{\alpha}r^{1-p}g(v_-^\eps
      +v_+^\eps),\  \ r>0,
\\
&v_\mp^\eps\big|_{t=0} = v_{0\mp}^\e\,,
\end{split}\right.
\end{equation}
where $ \operatorname{supp} v_{0\mp}^\e \subset [\delta^\e,+\infty[$
for some $\delta^\e >0$, and $v_{0\mp}^\e \in L^\infty(\R_+)$. Then by
finite speed of propagation, the term 
$r^{1-p}$ is harmless at least up to time $\delta^\e$. In that case,
local in time existence of solutions to \eqref{eq:pbreduit2} is easy. 
\begin{lem}
\label{lem:existsans}
Fix $\alpha \geq 0$, $p>1$. 
Let $v_{0\mp}^\e \in L^\infty(\R_+)$ such that $ \operatorname{supp}
v_{0\mp}^\e \subset [\delta^\e,+\infty[$ for some $\delta^\e >0$. Then
there exists $T^\e$, with $0<T^\e <\delta^\e$, and a unique solution
$(v_-^\e,v_+^\e) \in L^\infty\cap C([0,T^\e]\times\R_+)^2$ to the initial
value problem \eqref{eq:pbreduit2}. 
\end{lem}

When the incoming wave $v^\e_-$ reaches the origin $\{r=0\}$, a
boundary condition is needed in order to solve the above system. We
are interested in that given in \eqref{eq:pbreduit}, that is
$$(v_-^\eps +v^\eps_+)\big|_{r=0}=0.$$
Consider the mixed problem
\begin{equation}
\label{eq:pbreduit3}\left\{
  \begin{split}
      &(\d_t\pm \d_r)v_\pm^\eps = \eps^{\alpha}r^{1-p}g(v_-^\eps
      +v_+^\eps),\  \ r>0,
\\
 &  (v_-^\eps +v^\eps_+)\big|_{r=0}=0\,,\\
&v_\mp^\eps\big|_{t=0} = v_{0\mp}^\e.
\end{split}\right.
\end{equation}
The boundary condition 
compensates the 
singularity $r^{1-p}$ when $r$ gets close to zero. Indeed, Taylor's
formula yields, for $C^1$ solutions, 
$$(v_-^\eps +v^\eps_+)(t,r)= r\d_r (v_-^\eps +v^\eps_+)(t,r) +o(r),
\textrm{ as }r \to 0.$$
Now from the differential equation, we also have
$$\d_r (v_-^\eps +v^\eps_+) = \d_t (v_-^\eps -v^\eps_+),$$
so if we know that the time derivatives of $v^\e$ remain bounded, then
the singularity $r^{1-p}$ is compensated. We have precisely,
\begin{equation}\label{eq:poids}
|(v_-^\eps +v^\eps_+)(t,r)|\leq \frac{4r}{r+\e}\left(|v_-^\eps| +
 |v_+^\eps| +|\e \d_t v_-^\eps| +|\e \d_t v_+^\eps| \right)(t,r), \
\forall t,r\geq 0.
\end{equation}
This is the strategy we used
in \cite{CR2}, Proposition~3.4, to prove local existence, in the case
$p>2$. Notice that at this 
stage, the dependence upon $\e$ is unimportant, and that 
 $p>1$ suffices for local existence. 
\begin{lem}
\label{lem:localavec}
Fix $\alpha \geq 0$, $p>1$. 
Let $v_{0\mp}^\e \in W^{1,\infty}(\R_+)$ such that
\begin{equation}\label{eq:comp}
v_{0-}^\e(0)+v_{0+}^\e(0)=0\, ,\ \ \ \textrm{ and }\
\d_r v_{0-}^\e(0)-\d_rv_{0+}^\e(0)=0,
\end{equation}
then there exists $T^\e>0$ and a unique solution  $(v_-^\e,v_+^\e) \in
C^1\cap W^{1,\infty}([0,T^\e]\times\R_+)^2$ of \eqref{eq:pbreduit3}.
\end{lem}
As recalled in the beginning of this section, such a result will be
needed only in the proof of Theorem~\ref{th:sub}, and not in the
proof of Theorem~\ref{th:super}. From now on in this section, we
assume $\alpha >\max (0,p-2)$. In \cite{CR2}, we also proved
that the solutions of \eqref{eq:pbreduit3} with $\e=1$ and $p>2$ are
global provided that the initial data $v_{0\pm}$ are sufficiently small. 
\begin{prop}[\cite{CR2}, Proposition~3.5]
\label{lem:globalavec}
Fix $\alpha \geq 0$, $p>2$.  There are   constants $K_1$ and  $K_1^\prime>0$ 
so that for all
initial data $\psi_0\in C^1([0,\infty))$ satisfying
$$
\|\psi_0, \d_r\psi_0\|_{L^\infty([0,\infty[)} \le K_1\,,
$$
 and the compatibility conditions
\quad
$$
\psi_{0+}(0)+\psi_{0-}(0) \ = \ 0\,,
\qquad
{\rm and}
\qquad
\d_r\psi_{0+}(0)-\d_r \psi_{0-}(0) \ = \ 0\,,
$$
there is a unique solution $\psi\in C^1([-\infty,\infty]\times
[0,\infty[)$ of 
\begin{equation}
\label{eq:existpsi}
\left\{
\begin{split}
&(\d_t\pm \d_r)\psi_\pm\ = \ r^{1-p}g(\psi_-
+\psi_+), 
\\
&   \psi_-(t,0)+\psi_{+}(t,0)\ =\ 0
\,,
\\
&  \psi_{\mp}|_{t=0}\ =\ \psi_{0\mp}.
\end{split}
\right.
\end{equation}
In addition, 
\begin{equation*}
\|\psi, \partial_{t}\psi\|_{L^\infty([-\infty,\infty]\times[0,\infty[)}
 \le K^\prime_1\,
\|\psi_0, \d_r\psi_0\|_{L^\infty([0,\infty[)}
\,.
\end{equation*}
\end{prop}
The idea is then to find a scaling such that we can use the above
proposition to prove global existence for \eqref{eq:pbreduit} when
$\alpha>p-2>0$, along 
with useful estimates. Try
\begin{equation}\label{eq:scaling}
v^\e_\pm(t,r) =\e^\gamma \psi^\e_\pm(\tau,\rho)\big|_{\tau
={ \textstyle\frac{t}{\e}}, \rho={\textstyle\frac{r}{\e}}}\ \ , \ \ 
\psi^\e_\pm(\tau,\rho) = \e^{-\gamma}v^\e_\pm(\e \tau,\e \rho)    \ .
\end{equation}
Then $\psi^\e_\pm$ solve the differential equations,
$$(\d_\tau \pm \d_\rho)\psi^\e_\pm = \eps^{\alpha}(\e
\rho)^{1-p}\e^{(p-1)\gamma} g(\psi^\e_-+ \psi^\e_+).$$
Choose $\gamma$ so that the powers of $\e$ cancel,
$$\alpha +1-p+(p-1)\gamma =0 \Leftrightarrow \gamma =
1-\frac{\alpha}{p-1}\  .$$ 
Since $\gamma$ is negative, 
$$\|\psi^\e_\pm(0)\|_{L^\infty_\rho}\ ,\
\|\d_\tau\psi^\e_\pm(0)\|_{L^\infty_\rho} \Tend \e 0 0,$$
so we can apply Prop.~\ref{lem:globalavec}, for $\eps$ sufficiently
small. Moreover, it provides $L^\infty$ estimates for $v^\e_\pm, \e
\d_t v^\e_\pm$. We deduce that $\d_r v^\e_\pm\in L^\infty$ from the
differential equations and \eqref{eq:poids}.  
\begin{cor}\label{cor:globalp>2}
Assume $\alpha >p-2>0$. Then there exists $\e_0>0$ such that for
$0<\e\leq \e_0$, \eqref{eq:pbreduit} has a unique solution
$(v_-^\e,v_+^\e) \in
C^1\cap W^{1,\infty}(\R_+\times\R_+)^2$. Moreover, there exists $C$ such that
for any $\e\in ]0,\e_0]$, 
$$\|v^\e_\pm, \e \d_t v^\e_\pm\|_{L^\infty(\R_+\times\R_+)}\leq C.$$
\end{cor}
We now have to prove global existence when $1<p\leq 2$ and $\alpha
>0$. Since local existence is known (Lemma~\ref{lem:localavec}),
global existence is a consequence of \emph{a priori} estimates,
which follow from \eqref{eq:poids} and
Lemma~\ref{lem:linearest}. Define
$${\bf m}_\pm^\e(t)=\sup_{0\leq s\leq t}\left( 
\|v^\e_{\pm}(s)\|_{L^\infty}+ \|\e \d_t v^\e_{\pm}(s)\|_{L^\infty}
\right).$$
By assumption, ${\bf m}_\pm^\e(0)$ are bounded independent of $\e\in
]0,1]$. Lemma~\ref{lem:linearest} and \eqref{eq:poids} yield, for
$T>0$,
\begin{equation*}
\begin{split}
{\bf m}_-^\e(T) + {\bf m}_+^\e(T)& \leq  C\left( {\bf m}_-^\e(0) + {\bf
m}_+^\e(0) +
h_p(\e,T)\left( {\bf m}_-^\e(T) + {\bf m}_+^\e(T)
\right)^p\right)\\
& \leq  C_0+ C\left( 
h_p(\e,T)\left( {\bf m}_-^\e(T) + {\bf m}_+^\e(T)
\right)^p\right),
\end{split}
\end{equation*}
with 
\begin{equation*}
h_p(\e,T)=\left\{
\begin{split}
\e^{\alpha}T^{2-p}\ , & \textrm{ if }1<p<2,\\
\e^{\alpha}\log \left({1+{\textstyle{\frac{T}{\e}}}}\right)\ 
,& \textrm{ if }p=2.
\end{split}
\right.
\end{equation*}
In particular, for fixed $T>0$, $h_p(\e,T)\to 0$ as $\e \to
0$. Therefore, for any $T>0$, there exists $\e(T)>0$ such that for any
$\e\in ]0,\e(T)]$, $v^\e$ exists and, 
$${\bf m}_-^\e(T)+{\bf m}_+^\e(T)\leq 4C_0.$$
If ${\bf m}_-^\e(T)+{\bf m}_+^\e(T)\leq 4C_0$, it follows that
$$
{\bf m}_-^\e(T) + {\bf m}_+^\e(T)\leq C_0 +C h_p(\e,T) (4C_0)^p,
$$
and therefore for $0<\eps\leq\eps(T)$,
\begin{equation*}
{\bf m}_-^\e(T) + {\bf m}_+^\e(T)\leq 2 C_0
\end{equation*}
This proves that in fact ${\bf m}_-^\e(T)+{\bf m}_+^\e(T)< 4C_0 $
since if that were not true, 
there would be a first $T_*$ where ${\bf m}_-^\e(T)+{\bf m}_+^\e(T)=
4C_0$, and at 
that value of $T$ the above estimate leads to the contradiction
$4C_0<2C_0$.
\begin{prop}
\label{prop:globalp<=2}
Assume that $\alpha>0$ and $1<p\leq 2$. Let $T>0$. Then there exists $\e(T)>0$
such that for 
$0<\e\leq \e(T)$, \eqref{eq:pbreduit} has a unique solution
$$(v_-^\e,v_+^\e) \in
C^1\cap W^{1,\infty}([0,T]\times\R_+)^2.$$
Moreover, there exists $C$ such that
for any $\e\in ]0,\e(T)]$, 
$$\|v^\e_\pm, \e \d_t v^\e_\pm\|_{L^\infty([0,T]\times\R_+)}\leq C.$$
\end{prop}

%% file: sub.tex
\section{The subcritical case}
\label{sec:sub}

Now Theorem~\ref{th:sub} is a straightforward consequence of
Lemma~\ref{lem:linearest}, Corollary~\ref{cor:globalp>2} and  
Proposition~\ref{prop:globalp<=2}. In the statement of
Theorem~\ref{th:sub}, we distinguished three cases; $p>2$, $p=2$ and
$1<p<2$. The distinction $p{> \atop \leq}2$ appeared in the previous
section, in Corollary~\ref{cor:globalp>2} and
Proposition~\ref{prop:globalp<=2}. It corresponds to the question of the
integrability at infinity of the mapping $r\mapsto r^{1-p}$. The
further distinction $p=2$ corresponds to the local integrability of
this mapping. \\

Define the remainder $w^\e_\pm = v^\e_\pm - (v_\pm^\eps)_{\rm
free}$. It solves the mixed problem,
\begin{equation}
\label{eq:reste}
\left\{
  \begin{split}
      (\d_t\pm \d_r)w_\pm^\eps & = \eps^{\alpha}r^{1-p}g(v_-^\eps +v_+^\eps),
\\
(w_-^\eps +w^\eps_+)\big|_{r=0}& =0\,,\\
w_\mp^\eps\big|_{t=0} &= \pm \eps
P_1\left( r,\frac{r-r_0}{\eps}\right)\,,
\end{split}\right.
\end{equation}
Let $I$ be an interval of the form $[0,T[$, with $T\in \R_+\cup
\{+\infty\}$. From Lemma~\ref{lem:linearest}, we have
\begin{equation*}
\|w^\e_\pm\|_{L^\infty(I\times \R_+)}\leq 2\e\|P_1\|_{L^\infty} + 
2\sup_{t\in I} \sup_{\gamma \in \Gamma^t}\int_\gamma
\eps^{\alpha}r^{1-p}g(v_-^\eps +v_+^\eps).
\end{equation*}
Using \eqref{eq:poids}, we also have, 
\begin{equation}
\label{eq:sub}
\|w^\e_\pm\|_{L^\infty(I\times \R_+)}\leq C\e + 
C\left(\sup_{t\in I} \sup_{\gamma \in \Gamma^t}\int_\gamma
\frac{\eps^{\alpha}}{(r+\e)^{p-1}}\right)
\left\| v^\e,\e \d_t v^\e\right\|_{L^\infty(I\times \R_+)}^p.
\end{equation}
Differentiating \eqref{eq:reste} with respect to time yields, using
the differential equation \eqref{eq:reste} to find the initial data, 
\begin{equation*}
\left\{
  \begin{split}
      (\d_t\pm \d_r)\e\_d w_\pm^\eps & =
      \eps^{\alpha}r^{1-p}g'(v_-^\eps +v_+^\eps)\e\d_t(v_-^\eps +v_+^\eps), 
\\
\e\d_t(w_-^\eps +w^\eps_+)\big|_{r=0}& =0\,,\\
\e\d_tw_\mp^\eps\big|_{t=0} &= \eps
(\e\d_r P_1 +\d_z P_1)\left(
r,\frac{r-r_0}{\eps}\right)+\e^\alpha r^{1-p}g(v_-^\eps
+v_+^\eps)\big|_{t=0} \,,
\end{split}\right.
\end{equation*}
Since the initial data for $v^\e_\pm$ are supported in $|r-r_0|\leq
z_0\e$, the term $r^{1-p}$ in the initial data for $\e\d_tw_\mp^\eps$
is harmless regarding to $L^\infty$ estimates. 
From Lemma~\ref{lem:linearest} and \eqref{eq:poids}, we have, 
\begin{equation*}
\|\e\d_tw^\e_\pm\|_{L^\infty(I\times \R_+)}\lesssim \e^{\min(1,\alpha)} + 
\left(\sup_{t\in I} \sup_{\gamma \in \Gamma^t}\int_\gamma
\frac{\eps^{\alpha}}{(r+\e)^{p-1}}\right)
\left\| v^\e,\e \d_t v^\e\right\|_{L^\infty(I\times \R_+)}^p.
\end{equation*}

\noindent {\bf The case $p>2$.} Assume $\alpha >p-2>0$. Then we have
$$ \sup_{t\geq 0} \sup_{\gamma \in \Gamma^t}\int_\gamma
\frac{\eps^{\alpha}}{(r+\e)^{p-1}} \leq C\e^{\alpha +2-p},$$
and from Corollary~\ref{cor:globalp>2}, there exists $\e_0$ such that if
$\e \in ]0,\e_0]$, $\left\| v^\e,\e \d_t
v^\e\right\|_{L^\infty(\R_+\times \R_+)}$ is bounded. This proves the
first part of Theorem~\ref{th:sub}.\\

\noindent {\bf The case $p=2$.} When $\alpha >p-2=0$, then from
Prop.\ref{prop:globalp<=2}, for every $T>0$, there exists $\e(T)$ such
that if $\e\in ]0,\e(T)]$, $\left\| v^\e,\e \d_t
v^\e\right\|_{L^\infty([0,T[\times \R_+)}$ remains bounded. Moreover, for
any fixed $T>0$, we have
$$\sup_{\gamma \in \Gamma^T}\int_\gamma
\frac{\eps^{\alpha}}{r+\e} \leq C \eps^{\alpha} \log\frac{T}{\e}\, .$$

\noindent {\bf The case $1<p<2$.} When $\alpha>0$ and $1<p<2$, then
the only difference with the previous case is that the mapping
$r\mapsto r^{1-p}$ is locally integrable, hence the bound
$$\sup_{\gamma \in \Gamma^T}\int_\gamma
\frac{\eps^{\alpha}}{(r+\e)^{p-1}}\leq \sup_{\gamma \in \Gamma^T}\int_\gamma
\frac{\eps^{\alpha}}{r^{p-1}}
 \leq C \eps^{\alpha} T^{2-p}\, .$$
This completes the proof of Theorem~\ref{th:sub}.\qed

%% file: super.tex
\section{The supercritical case}
\label{sec:super}
We conclude by proving Theorem~\ref{th:super}, using
Proposition~\ref{prop:approx}. 

Before going into details, we explain how to
construct the approximate solutions that  lead to the
result. Assume for instance that $\alpha >0$ and, since we are in a
supercritical case, $p-2>\alpha$. The hypothesis $\alpha>0$ suggests
that the nonlinear term in \eqref{eq:pbreduit} is negligible when $r$
is not too small. This is the outline of the proof of
Theorem~\ref{th:sub}, and we could prove this way that $(v_\pm^\eps)_{\rm
free}$ is a good approximation for $v^\e_\pm$ at least for $r_0-t \gg
\e^{\frac{\alpha}{p-2}}$. This boundary layer is larger than in the
critical case $\alpha=p-2>0$ studied in \cite{CR2}, and nonlinear
effects possibly occur sooner. Considering the case $\alpha=0$ gives
us a further hint. We proved in \cite{CR} that a good approximation at
leading order was given by
$$(\d_t\pm \d_r)(v_\pm^\eps)_{\rm app}=r^{1-p}g\left((v_\pm^\eps)_{\rm
app}\right), $$ 
that is, solutions of ordinary differential equations along the rays
of geometrical optics. A natural generalization of this approach to
the case $\alpha \geq 0$ leads to 
\begin{equation}
\label{eq:app}
\left\{
  \begin{split}
      &(\d_t\pm \d_r)(v_\pm^\eps)_{\rm app} =
      \eps^\alpha r^{1-p}g((v_\pm^\eps)_{\rm app}), \ \ (t,r)\in \left([0,
      r_0-z_0\e[\times \R^*_+\right)\, ,\\
 &  (v_\pm^\eps)_{\rm app}\big|_{t=0} = P_\pm\left
   ( r,\frac{r-r_0}{\eps}\right).
   \end{split}
\right.
\end{equation}
We consider the region $t< r_0-z_0\e$ so that no boundary condition
is needed on $\{r=0\}$, since the compact support of $P_\pm$ and the
finite speed of propagation make $(v_\pm^\eps)_{\rm app}$ zero in the
region $\{ r+t < r_0-z_0\e\}$ (see Figure~\ref{fig:ogg}). 

Introduce the function $F_p$, defined for $y\geq x > 0$ by,
\begin{equation*}
F_p(x,y)=\int_x^y \frac{ds}{s^{p-1}}=\left\{
\begin{split}
&\frac{1}{p-2}(x^{2-p}-y^{2-p}), \textrm{ if } p>2,\\
&\log \frac{y}{x}, \textrm{ if } p=2.
\end{split}\right.
\end{equation*}
Then \eqref{eq:app} can be solved explicitly, 
\begin{equation*}
\begin{split}
(v_-^\eps)_{\rm app}(t,r) &= \frac{P_-(r+t,z)}{\Bigl( 
1+a2^{-p}(p-1)\eps^\alpha F_p(r,r+t)|P_-(r+t,z)|^{p-1}\Bigr)^{\frac{1}{p-1}}}
\Biggl|_{\textstyle z=\frac{r+t-r_0}{\eps}} ,\\
(v_+^\eps)_{\rm app}(t,r) &=\frac{P_+(r-t,z)}{\Bigl( 
1+a2^{-p}(p-1)\eps^\alpha F_p(r-t,r)|P_+(r-t,z)|^{p-1}\Bigr)^{\frac{1}{p-1}}}
\Biggl|_{\textstyle z=\frac{r-t-r_0}{\eps}}.
\end{split}
\end{equation*}
We have, 
$$\operatorname{supp}(v_-^\eps)_{\rm app}=\{ |r+t-r_0|\leq
z_0\eps\},\ 
\operatorname{supp}(v_+^\eps)_{\rm app}=\{ |r-t-r_0|\leq
z_0\eps\}.$$
Since for any $y>0$ and any $p\geq 2$, $F_p(x,y)\Tend x 0
{+\infty}$, the sign of $a$ is crucial. 
\begin{itemize}
\item If $a>0$ (dissipative case), then 
$(v_-^\eps)_{\rm app}$ tends to zero before reaching the focus. 
\item If $a<0$ (accretive case), then both $(v_-^\eps)_{\rm
app}$ and $(v_+^\eps)_{\rm app}$ blow up in finite time. 
\end{itemize}
In the accretive case, $(v_+^\eps)_{\rm app}$ may blow
up at time $T^*<r_0$, that is, before any focusing. To avoid useless
distinctions, we assume $P_+\equiv 0$, so that the only
``interesting'' phenomena occur when $t$ approaches $r_0$. 
 
The explicit formulae show that when $p>2$, for
$\e^\alpha r^{2-p}\gg 1$, 
$$(v_\pm^\eps)_{\rm app}(t,r)=(v_\pm^\eps)_{\rm free}(t,r)+O(\e^\alpha
r^{2-p}).$$
Assume first $p-2>\alpha \geq 0$. Let $\lambda >0$. We prove that
$(v_\pm^\eps)_{\rm app}$ gives a 
good approximation of the exact solution in the region $r\geq \lambda
\eps^{\frac{\alpha}{p-2}}$ before the focus. 
Notice that in this
region, $(v_\pm^\eps)_{\rm app}$ and $(v_\pm^\eps)_{\rm free}$ have
ceased to be close to each other; 
the nonlinear effects 
are significant.  By finite
propagation speed, this area is defined by 
$$t\leq r_0-\lambda \eps^{\frac{\alpha}{p-2}}-z_0\eps =
T_{\lambda,\eps}.$$ 
We prove that $(v_\pm^\eps)_{\rm app}$ remains a good
approximation of the exact solution up to time $T_{\lambda,\eps}$. 
\begin{prop}
\label{prop:solapp2}
Assume $p-2>\alpha \geq 0$, and let $\lambda >0$. Then
$(v_\pm^\eps)_{\rm app}$ is a  
good approximation of the exact solution at least for $t\in [0,
T_{\lambda,\eps}]$,
$$\|v^\eps_\pm - (v_\pm^\eps)_{\rm app}\|_{L^\infty([0,
T_{\lambda,\eps}]\times\R_+)}=O\left(
\eps^{1-\frac{\alpha}{p-2}}\right)\ \  \Big(O(\e)\ \textrm{ if }p=2\Big).$$ 
\end{prop}
\begin{proof}[Proof of Proposition~\ref{prop:solapp2}]
The proof relies on Proposition~\ref{prop:approx}, which we apply as in
Example~\ref{ex:gener}. 

By finite propagation speed, both the exact and approximate solution
are zero in the region 
$$\Big\{(t,r)\in [0,T_{\lambda,\eps}]\times \R_+\, ;\, \, r+t\leq
T_{\lambda,\eps}\Big\}.$$
Define the remainder $w_\pm^\eps =v^\eps_\pm- (v_\pm^\eps)_{\rm
app}$. Before rays reach the boundary $\{r=0\}$, 
\begin{equation*}
\left\{
  \begin{split}
      (\d_t\pm \d_r)w_\pm^\eps &=
      \eps^\alpha r^{1-p}\left( g(v_-^\eps +v_+^\eps)-
      g((v_\pm^\eps)_{\rm app})\right), \\ 
 w_\pm^\eps\big|_{t=0} &=\mp\eps  P_1\left( r,\frac{r-r_0}{\eps}\right).
   \end{split}
\right.
\end{equation*}
In order to apply Proposition~\ref{prop:approx}, write the right hand
side as 
\begin{equation*}
\begin{split}
\eps^\alpha r^{1-p}\left(g(v_-^\eps +v_+^\eps)-
      g((v_\pm^\eps)_{\rm app})\right) =& \eps^\alpha r^{1-p} \Big(
g(v_-^\eps +v_+^\eps) - 
g((v_-^\eps)_{\rm app}+(v_+^\eps)_{\rm app} )\\
&+g((v_-^\eps)_{\rm app}+(v_+^\eps)_{\rm app} )
      -g((v_\pm^\eps)_{\rm app})\Big).
\end{split}
\end{equation*}
Using Taylor's Theorem, the first term satisfies
$$\eps^\alpha r^{1-p}\left(g(v_-^\eps +v_+^\eps) - 
g((v_-^\eps)_{\rm app}+(v_+^\eps)_{\rm app} )\right) =
\eps^\alpha r^{1-p}(w^\eps_-+w^\eps_+) f^\eps(t,r),$$
where $f^\eps$ is uniformly bounded on any set on which the families
$v^\eps$ and $ (v^\eps)_{\rm app}$ are uniformly bounded. The point
now is that this term has exactly the properties  mentioned in
Example~\ref{ex:gener}; in particular, the assumptions of
Proposition~\ref{prop:approx} are satisfied up to time
$T=T_{\lambda,\e}$. 

Define the source term 
$$S_\pm^\eps (t,r) \ := \ \eps^\alpha r^{1-p} \left( g((v_-^\eps)_{\rm
app}+(v_+^\eps)_{\rm app} ) 
      -g((v_\pm^\eps)_{\rm app})\right).$$
Then $\operatorname{supp}S_\pm^\eps \subset
\operatorname{supp}(v_\mp^\eps)_{\rm app}$ , and
\begin{equation}
\left\{
  \begin{split}
      (\d_t\pm \d_r)w_\pm^\eps &=
      \eps^\alpha r^{1-p}f^\eps(t,r)(w^\eps_-+w^\eps_+) + S_\pm^\eps, \\ 
 w_\pm^\eps\big|_{t=0} &=\mp\eps  P_1\left
   ( r,\frac{r-r_0}{\eps}\right).
   \end{split}
\right.
\end{equation}
From the explicit expression of $(v^\eps_\pm)_{\rm app}$, we see that
for any fixed $\lambda >0$,
there exists $C_\lambda>0$ independent of $\eps$ such that
\begin{equation}
\big\|(v^\eps_\pm)_{\rm app}
\big\|_{L^\infty([0,T_{\lambda,\eps}]\times\R_+)}\ \leq \ C_\lambda.
\end{equation}
We shall prove that
there is $\eps(\lambda)>0$, so that for
$\eps\in ]0,\eps(\lambda)]$, 
\begin{equation}
\label{tiger}
\|w^\eps_\pm\|_{L^\infty([0,T_{\lambda,\eps}]\times\R_+ )}
\ \le \ C_\lambda\,\eps^{1-\frac{\alpha}{p-2}}\, \qquad
\Big(C_\lambda\,\eps\ \textrm{ if }p=2\Big)\, .
\end{equation}
This implies the  
error estimate of
Proposition~\ref{prop:solapp2}. For the sake of readability, we will
omit the distinction $p=2$, and keep the notation
$\eps^{1-\frac{\alpha}{p-2}}$, with an obvious convention. 

From Lemma~\ref{lem:existsans}, $v^\e$, hence $w^\e$, is defined,
bounded and continuous, locally in time. At time $t=0$, it is of order
$\e$, so there exists $t^\e>0$ such that
$$\|w^\eps_\pm\|_{L^\infty([0,t^\eps]\times\R_+ )}< 2\e
\|P_1\|_{L^\infty},$$
and we have, possibly increasing the
value of $C_\lambda$, 
$$\|w^\eps_\pm\|_{L^\infty([0,t^\e]\times\R_+ )}
\ \le \ C_\lambda\,\eps^{1-\frac{\alpha}{p-2}}\,.$$
So long as $w^\eps_\pm$ is pointwise bounded by $2C_\lambda$,
$f^\eps$ remains uniformly bounded. As noticed in
Example~\ref{ex:gener}, there exists $C_2(\lambda)$ 
such that 
$$\int_0^{T_{\lambda,\eps}}\sup_{\gamma_\in
\Gamma^t}
\eps^\alpha r^{1-p}\left| f^\eps\right| dt\ \leq \ 
C_2(\lambda). $$
Now consider the source terms. On $\operatorname{supp}
(v_+^\eps)_{\rm app}$ (of size $2z_0\eps$, and transverse to any $\gamma
\in \Gamma_-^{T(\lambda,\e)}$), $r>r_0/2$, therefore the
singular term $r^{1-p}$ is bounded and
$$\sup_{\gamma_- \in \Gamma_-^{T(\lambda,\e)}}\int_{\gamma_-}
|S^\eps_-|\leq C\eps^{1+\alpha}.$$ 
More delicate is the treatment of $S^\eps_+$. Because $\operatorname{supp}
(v_-^\eps)_{\rm app}$ is of size $2z_0\eps$ and transverse to any $\gamma_+
\in \Gamma_+^{T(\lambda,\e)}$, one has 
$$\sup_{\gamma_+ \in \Gamma_+^{T(\lambda,\e)}}\int_{\gamma_+} |S^\eps_+|\leq
C \eps^\alpha \int_{\lambda \eps^{\frac{\alpha}{p-2}}}^{\lambda
\eps^{\frac{\alpha}{p-2}} +2z_0\eps} r^{1-p}dr= C \eps^\alpha F_p\left(\lambda
\eps^{\frac{\alpha}{p-2}}, \lambda
\eps^{\frac{\alpha}{p-2}} +2z_0\eps \right) .$$ 
If $p>2$,
\begin{equation*}
\begin{split}
(p-2)\eps^\alpha F_p\left(\lambda
\eps^{\frac{\alpha}{p-2}} , \lambda
\eps^{\frac{\alpha}{p-2}}+2z_0\eps \right)=& \eps^\alpha\left(\left(\lambda
\eps^{\frac{\alpha}{p-2}} \right)^{2-p} - \left(\lambda
\eps^{\frac{\alpha}{p-2}}+2z_0\eps\right)^{2-p}\right)\\
=& \lambda^{2-p} \left(1-\left(1
+\frac{2z_0}{\lambda}\eps^{1-\frac{\alpha}{p-2} } 
\right)^{2-p}\right)\\
=& O\left(\lambda^{1-p} \eps^{1-\frac{\alpha}{p-2} }\right).
\end{split}
\end{equation*}
The case $p=2$ yields the same estimate.
From Proposition~\ref{prop:approx}, $w^\eps_\pm$ remains pointwise
bounded by $2C_\lambda$ for $t\in [0,T(\lambda,\e)]$, for $\e \in
]0,\e(\lambda)]$. This yields \eqref{tiger}, and completes the proof of
Proposition~\ref{prop:solapp2}. 
\end{proof}

Proposition~\ref{prop:solapp2} describes the behavior of
the exact solution $v^\e$ up to time $T(\lambda,\e)$. 
The outline of the end of the proof in the dissipative case is as
follows. At time $t=T(\lambda,\e)$, the approximate solution
$(v^\eps)_{\rm app}$ is of order $\lambda^{(p-2)/(p-1)}$ if
$p>2$ ($1/|\log \lambda|$ if $p=2$). Letting  $\e$ go to zero, with
$\lambda$  sufficiently small, shows that $v^\e$ becomes arbitrarily
small when approaching the focal point. Since the equation is
dissipative, this means that $v^\e$ is absorbed. The end of the proof
of Theorem~\ref{th:super} relies on energy estimates. 
For $q\geq 1$, define
$$
g_{q-1}(s) 
\ :=\ 
{d\over ds}\, |s|^{q}
\ = \ 
q|s|^{q-1}\,{\rm sgn}\,s\,.
$$
Then $g_{q}$ is a non-increasing odd function of $s$ which is homogeneous
of degree $q$.

In \eqref{eq:pbreduit}, multiply the equation satisfied by $v^\e_-$
by
$g_{q-1}(v^\e_-)$, 
and  the equation satisfied by $v^\e_+$ by 
$g_{q-1}(v^\e_+)$. 
Summing yields
\begin{equation*}
\begin{split}
\d_t (|v^\e_-|^q+|v^\e_+|^q)& +\d_r (|v^\e_+|^q-|v^\e_-|^q) = \\
& -a2^{-p} r^{1-p}\,
\big(g_{q-1}(v^\e_-) + g_{q-1}(v^\e_+)\big)\,
g_p(v^\e_-+v^\e_+)/(p+1)\,.
\end{split}
\end{equation*}
The signs of both $g_{q-1}(v^\e_-) + g_{q-1}(v^\e_+)$ and 
$g_p(v^\e_-+v^\e_+)$ are equal to the sign of the larger
of $v^\e_\pm$.  
Therefore, when $a>0$ (dissipative case), 
\begin{equation}
\label{eq:nrj}
\d_t (|v^\e_-|^q + |v^\e_+|^q)
\ -\ 
\d_r
(|v^\e_-|^q - |v^\e_+|^q)
\ \leq\ 
 0.
\end{equation}
For a fixed $t>0$, integrate this inequality from $r=0$ to
$r=t$. Recall that when $r=0$, we have $|v^\eps_-|=|v^\eps_+|$, so
this yields
$$\int_0^t  \d_t (|v^\eps_-|^q + |v^\eps_+|^q)(t,r)dr -
|v^\eps_-|^q(t,t) +|v^\eps_+|^q(t,t)\leq 0. $$
Therefore,
\begin{equation}\label{eq:zxc}
\d_t \int_0^t  (|v^\eps_-|^q + |v^\eps_+|^q)(t,r)dr\leq
2|v^\eps_-|^q(t,t).
\end{equation} 
Let $\lambda >0$. From Proposition~\ref{prop:solapp2}, for $0\leq
t\leq T_{\lambda,\eps}$,  
$$|v^\eps_-|^q(t,t)=|(v^\eps_-)_{\rm app}|^q(t,t)+O\left
( \eps^{q\left(1-\frac{\alpha}{p-2}\right)} \right).$$
For $t>\frac{r_0+z_0\eps}{2}$, $(v^\eps_-)_{\rm app}(t,t)=0$ and 
$$|v^\eps_-|^q(t,t)=O\left
( \eps^{q\left(1-\frac{\alpha}{p-2}\right)} \right).$$
Let $T>T_{\lambda,\eps}$ and integrate \eqref{eq:zxc} between $t=
T_{\lambda,\eps} $ and $t=T$. 
\begin{equation*}
\begin{split}
\|v_-^\eps(T)\|_{L^q(0,T)}^q +\|v_+^\eps(T)\|_{L^q(0,T)}^q \leq & 
 \|v_-^\eps(T_{\lambda,\eps})\|_{L^q(0,T_{\lambda,\eps})}^q +
\|v_+^\eps(T_{\lambda,\eps})\|_{L^q(0,T_{\lambda,\eps})}^q \\
& +C(\lambda)^q T\eps^{q\left(1-\frac{\alpha}{p-2}\right)},
\end{split}
\end{equation*}
and 
\begin{equation*}
\begin{split}
\|v_-^\eps(T)\|_{L^q(0,T)} +\|v_+^\eps(T)\|_{L^q(0,T)} \leq &
 \|v_-^\eps(T_{\lambda,\eps})\|_{L^q(0,T_{\lambda,\eps})} +
\|v_+^\eps(T_{\lambda,\eps})\|_{L^q(0,T_{\lambda,\eps})} \\
& +
C(\lambda)T^{1/q} \eps^{1-\frac{\alpha}{p-2}}.
\end{split}
\end{equation*}
Letting $q\rightarrow \infty$ yields
\begin{equation*}
\begin{split}
\|v_-^\eps(T)\|_{L^\infty(0,T)} +\|v_+^\eps(T)\|_{L^\infty(0,T)} \leq &
 \|v_-^\eps(T_{\lambda,\eps})\|_{L^\infty(0,T_{\lambda,\eps})} +
\|v_+^\eps(T_{\lambda,\eps})\|_{L^\infty(0,T_{\lambda,\eps})} \\
& +C(\lambda) \eps^{1-\frac{\alpha}{p-2}}.
\end{split}
\end{equation*}
Using Proposition~\ref{prop:solapp2} again, and the fact that
$\operatorname{supp}(v^\eps_+)_{\rm app}\subset \{|r-t-r_0|\leq
z_0\eps\}$, we also have
$$ \|v_-^\eps(T)\|_{L^\infty(0,T)} +\|v_+^\eps(T)\|_{L^\infty(0,T)}
\leq \|(v^\eps_-)_{\rm
app}(T_{\lambda,\eps})\|_{L^\infty(0,T_{\lambda,\eps})} +C(\lambda)
\eps^{1-\frac{\alpha}{p-2}}.$$ 
From the explicit expression for $(v^\eps_-)_{\rm app}$, 
\begin{equation}
\|(v^\eps_-)_{\rm
app}(T_{\lambda,\eps})\|_{L^\infty(0,T_{\lambda,\eps})}\leq
 C\times\left\{ 
 \begin{split}
  \lambda^{\frac{p-2}{p-1}} \ ,&\ \textrm{ if }p>2,\\
  1/|\log \lambda|\ ,&\ \textrm{ if }p=2,
 \end{split}\right.
\end{equation}
where $C$ does not depend on $\lambda$. Thus if $p>2$, 
$$ \|v_-^\eps(T)\|_{L^\infty(0,T)} +\|v_+^\eps(T)\|_{L^\infty(0,T)}
\leq C \lambda^{\frac{p-2}{p-1}} +C(\lambda)
\eps^{1-\frac{\alpha}{p-2}}.$$ 
Therefore,  for any $T\geq r_0$, and any $\lambda >0$,
\begin{equation}
\limsup_{\eps \rightarrow 0}\left(\|v_-^\eps(T)\|_{L^\infty(0,T)}
+\|v_+^\eps(T)\|_{L^\infty(0,T)}\right) \leq C \lambda^{\frac{p-2}{p-1}}.
\end{equation}
Letting $\lambda \rightarrow 0$ yields the first part of
Theorem~\ref{th:super} for $p>2$. The case $p=2$ is
straightforward. \\

When the equation is accretive ($a<0$), the energy estimate \eqref{eq:nrj}
becomes
\begin{equation}\label{eq:nrj6b}
\d_t (|v^\eps_-|^q + |v^\eps_+|^q)(t,r) -\d_r (|v^\eps_-|^q -
|v^\eps_+|^q)(t,r)\geq 0.
\end{equation}
Since we assumed $P_+\equiv 0$, 
$(v_-^\eps)_{\rm app}$ blows up in finite time, while
$(v_+^\eps)_{\rm app}$ does not.
The mechanism occurs quite the same way as the cancellation of
$(v_-^\eps)_{\rm app}$  in the dissipative case. For a fixed $t>0$,
integrating \eqref{eq:nrj6b} between $r=0$ and $r=\infty$ yields, by
finite speed of propagation, 
$$\int_0^\infty  \d_t (|v^\eps_-|^q + |v^\eps_+|^q)(t,r)dr=\d_t
\int_0^\infty   (|v^\eps_-|^q + |v^\eps_+|^q)(t,r)dr 
\geq 0. $$
Let $r_0>T>T_{\lambda,\eps}$ and integrate the above inequality between $t=
T_{\lambda,\eps} $ and $t=T$. 
\begin{equation*}
\|v_-^\eps(T)\|_{L^q(\R_+)}^q +\|v_+^\eps(T)\|_{L^q(\R_+)}^q \geq 
 \|v_-^\eps(T_{\lambda,\eps})\|_{L^q(\R_+)}^q +
\|v_+^\eps(T_{\lambda,\eps})\|_{L^q(\R_+)}^q .
\end{equation*}
Letting $q \to \infty$ yields,
\begin{equation*}
\|v_-^\eps(T)\|_{L^\infty(\R_+)} +\|v_+^\eps(T)\|_{L^\infty(\R_+)} \geq 
 \|v_-^\eps(T_{\lambda,\eps})\|_{L^\infty(\R_+)}.
\end{equation*}
From Proposition~\ref{prop:solapp2}, it follows
\begin{equation*}
\liminf_{\e \to 0}\left(\|v_-^\eps(T)\|_{L^\infty(\R_+)} +
\|v_+^\eps(T)\|_{L^\infty(\R_+)}\right) \geq  
 \liminf_{\e \to 0}\|(v_-^\e)_{\rm
 app}(T_{\lambda,\eps})\|_{L^\infty(\R_+)}.
\end{equation*}
Letting $\lambda \to 0$ yields the last part of
Theorem~\ref{th:super} with $T^*=r_0$, using the explicit form of
$(v_-^\e)_{\rm app}$. As we already mentioned, the result may hold
with $T^*<r_0$ if $(v_+^\e)_{\rm app}$ blows up before $(v_-^\e)_{\rm
app}$, in which case the proof is essentially as above.